\documentclass[letterpaper, 10 pt, conference]{ieeeconf}  
\usepackage{header}

\IEEEoverridecommandlockouts                              

\title{\LARGE \bf
Optimal, centralized dynamic curbside parking space zoning 
}

\author{Nawaf Nazir$^{1}$, Shushman Choudhury$^{2}$, Stephen Zoepf$^{2}$, Ke Ma$^{3}$, and Chase Dowling$^{1}$

\thanks{$^{1}$Nawaf Nazir and Chase Dowling are with Pacific Northwest National Laboratory \{nawaf.nazir,chase.dowling\}@pnnl.gov}%
\thanks{$^{3}$Ke Ma is with ISO New England }%
\thanks{$^{2}$Shushman Choudhury and Stephen Zoepf are with Lacuna Technologies. \{shushman.choudhury, stephen.zoepf\}@lacuna.ai}
}

\begin{document}

\maketitle
\thispagestyle{empty}
\pagestyle{empty}

\begin{abstract}

In this paper we formulate a dynamic mixed integer program for optimally zoning curbside parking spaces subject to transportation policy-inspired constraints and regularization terms. First, we illustrate how given some objective of curb zoning valuation as a function of zone type (e.g., paid parking or bus stop), dynamically rezoning involves unrolling this optimization program over a fixed time horizon. Second, we implement two different solution methods that optimize for a given curb zoning value function. In the first method, we solve long horizon dynamic zoning problems via approximate dynamic programming. In the second method, we employ Dantzig-Wolfe decomposition to break-up the mixed-integer program into a master problem and several sub-problems that we solve in parallel; this decomposition accelerates the MIP solver considerably. We present simulation results and comparisons of the different employed techniques on vehicle arrival-rate data obtained for a neighborhood in downtown Seattle, Washington, USA.

\end{abstract}

\section{Introduction}\label{sec:introduction}

The proliferation of smart-cities sensing devices---ranging from sensors mounted on civil infrastructure to crowd-sourced mobile data---has resulted in a wealth of new insights into how people interact with the built environment \cite{abel2021curbside}. In particular, data on curb activities like parking, passenger pickup-dropoff, and courier services are being collected by a variety of sensor types, such as cameras \cite{blakeley2013time}, proximity sensors \cite{ranjbari2021testing}, or indirectly via digitally logged transactions\cite{saha2019project}. These data have been used primarily to answer questions about demand for curb space~\cite{fiez2018data}; this work takes a forward looking perspective. This data will be useful as uses for curb space become increasingly varied, such as electric vehicle charging and micromobility services. For cities that will continue to centrally manage their curb zoning rules to adapt to an increasingly dynamic curb environment, we propose utilizing this proliferation of sensors to enable increasingly dynamic and flexible curb zoning.

Dynamic curb zoning implies changing zoning types, limitations and prices over short time horizons. We focus on daily time horizons; extant examples include curb zones available for parking during off-peak travel hours, and closed during peak travel hours to allow the curb lane to service more traffic~\cite{roe2017curb}. Other examples include airport departure zones serving arrivals and vice versa \cite{harris2017mesoscopic} and designating times for commercial loading \cite{giron2019commercial}. Typically these zone changes are reflective of periods of peak demand for varying use types, mitigating congestion at the curb or on the roadway~\cite{yu2021management}.

In this work we propose a framework for optimal, dynamic curbside zoning given information about time-varying curb demand. Given increasingly available sensors, we can better measure demand by time of day and vehicle type and thus rezone curbside space dynamically with respect to that demand. In this paper we formulate this dynamic rezoning problem as a mixed-integer program (MIP) that takes into account various transportation policy constraints derived from various sources of spatio-temporal data~\cite{fiez2018data}. To deal with the computational intractability that may arise due to the large number of integer variables in such problems, we investigate two different techniques to solve the MIPs in a scalable manner.

In the first approach, we use approximate dynamic programming (ADP) to solve the dynamic curb zoning problem. Dynamic programming has been applied successfully in many areas to solve mixed-integer programs~\cite{schwaeppe2019equal}. However, dynamic programming methods suffer from the curse of dimensionality and are not scalable with the increase in system size. Recent works in literature have looked into various approximate dynamic programming techniques that can improve scalability and computation time~\cite{lewis2009reinforcement,powell2007approximate}. In this work, we apply different ADP techniques to the dynamic curb zoning problem and compare the different techniques in terms of accuracy and computation time.

In the second approach, we solve the MIP with Dantzig-Wolfe decomposition, which is a delayed column generation technique~\cite{barnhart1998branch}. Several decomposition methods exist in literature to solve MIPs~\cite{galati2010decomposition}, however, the nature of the curb zoning problem allows us to break up the problem into a master problem and a set of sub-problems that can be solved in parallel. This problem structure allows the Dantzig-Wolfe decomposition to obtain a solution that scales well with problem size~\cite{andrianesis2021computation}. We baseline the two approaches against the global optimum obtained by solving the original MIP, and provide recommendations on how to choose the appropriate approach.

The rest of this paper is organized as follows: in Section~\ref{sec:background} we summarize existing research and describe some examples of dynamic curb space zoning. Section~\ref{sec:methods} formally defines the problem as a time-dynamic optimization program for sequences of zonings and solves the problem as a mixed-integer program. Further, we construct an example objective function for valuing individual curb spaces based on a hypothetical point-of-interest distribution and existing vehicle arrival rate data. In Section~\ref{sec:ADP} we formulate the curb zoning problem through approximate dynamic programming and compare different ADP techniques with the MIP solution. Section~\ref{sec:DW} employs the Dantzig-Wolfe decomposition method to solve the curb zoning problem and compare the results with the MIP solution. In Section~\ref{sec:Discussion} we discuss these different methods with respect to technological shifts in smart cities currently on the horizon, illustrating how this problem formulation is a valuable prospective avenue in which to view curb zoning as a land-use problem that need not be statically constrained. Finally, Section~\ref{sec:Conclusions} concludes our paper and lays out the scope for future work.

\section{Background}\label{sec:background}

Curb space is designated for users by various transit modes with categorical labels: e.g., bus stops, paid parking, commercial vehicle loading zones (CVLZ), passenger pick-up/drop-off  zones, curbside electric vehicle (EV) charging, and so on. Each usage has value, but one-to-one comparisons are exceptionally difficult, thus cities have used historical experience to balance varying objectives like access, congestion reduction, and productivity when \textbf{allocating} individual curb spaces to one or more of these modalities~\cite{guide2019introducing}. With higher resolution data streams and a better understanding of user elasticities~\cite{ostermeijer2021citywide}, adapting these allocations to changing demand in time and usage is becoming increasingly feasible. 

The literature on the subject of parking space assignment and curb zoning can be broken down into two fundamental perspectives: centralized and decentralized---from the viewpoint of the municipality that owns the curb real estate.

\begin{enumerate}
    \item \textbf{Centralized}: a municipality or curb manager maintains total control over zoning and space allocation and must solve for a network-wide optimal set of allocations \cite{inci2015review}.
    \item \textbf{Decentralized}: a market or auctioning mechanism is utilized that allows competing actors to bid for curb space access with varying degrees of oversight/constraints (i.e. bidding between modalities---people parking competing with trucks needing CVLZ access, or bidding amongst a single modality---commercial delivery companies bidding against one another for rights to access a space)~\cite{kong2018iot}.
\end{enumerate}

In this paper we focus on the centralized case of curb zoning. This is distinct from centralized space assignment where individual vehicles are routed to determined locations \cite{nakazato2019parking}; rather, the zoning problem seeks to determine where and when these parking locations should be made available and to what vehicle types. 

We can represent our problem as an instance of general multi-agent planning with shared resources~\cite{torreno2017cooperative}. The standard multi-agent planning problem is typically formulated as a combinatorial optimization or MIP, factored over the agents. Each agent has its local constraints and objective function (the global objective is simply the sum of agent objectives), and agents interact with each other through shared resources that act as multi-agent constraints. Prior research presented an optimal multi-agent planning algorithm based on Dantzig-Wolfe decomposition and Gomory cuts~\cite{hong2011optimal}. The second of our proposed approaches builds upon the essence of their work, albeit for solving our problem that has more difficult constraints.

\section{Methods}\label{sec:methods}


We now discuss how to solve for an optimal sequence of zonings under various constraints derived from transportation-oriented policies, assuming that we have an appropriate way to value curb spaces over time. We then construct an example curb valuation function (based on real-world arrival rate data) that serves as an objective function for our simulated results.

\subsection{MIP Formulation}

Denote $u_{k,c_j}^{(i)}\in \mathcal{U}$ as the curb allocation variable (decision variable) at time-step $k$, for curb $c_j$ and curb-type $(i)$, where $k\in \mathcal{T}= \{1,\hdots, T\}$, $c_j \in \mathcal{N}= \{c_1,\hdots, c_N\}$ and $(i) \in \mathcal{M}=\{\text{pp,cv,bus}\}$, with `pp' representing paid parking allocation, `cv' representing commercial vehicle allocation and `bus' representing allocation for public transportation. $u_{k,c_j}^{(i)}$ is a binary variable with $u_{k,c_j}^{(i)}=1$ if allocated and $u_{k,c_j}^{(i)}=0$ if not allocated. We can extend our formulation to more curb allocation types without loss of generality. 

We want to maximize the sum of finite known sequences of bounded functions $F_k:\mathcal{U}\rightarrow \mathbb{R}$, describing the utility of the curb space as a function of the allocation. In our case we define $F_k$ as a linear function of the curb allocation variable $u_{k,c_j}^{(i)}$ as:
\begin{align}
    F_k(u_{k,c_j},c_j)=\sum_{i \in \mathcal{M}}H_{k,c_j}^{(i)}u_{k,c_j}^{(i)}\label{eq:obj_def}
\end{align}
where $H_{k,c_j}^{(i)}$ is the normalized valuation for curb $c_j$ at time-step $k$ and for curb allocation type $i$. We assume that we know this value beforehand because we have calculated it from the available data on curb usage.

Apart from maximizing this function, we are also required to meet certain constraints with respect to the curb allocation. Some of them are listed below. 

\begin{enumerate}
    \item Curb changes: We would like to bound the number of curb allocation changes from one time-step to the next. Enforcing this constraint avoids unnecessary switching in curb space allocations over time. It can be expressed as:
    \begin{align}
    \sum_{c_j \in \mathcal{N}}\sum_{i\in \mathcal{M}}\frac{1}{2}\lVert u_{k+1,c_j}^{(i)}-u_{k,c_j}^{(i)} \rVert_1 \le b \ \forall k \in \mathcal{T}\backslash T \label{eq:change_cons}
    \end{align}
    \item Number of allocations at a time: We impose upper and lower bounds on the number of curb allocation types at any particular time, e.g., we would like a minimum and maximum number of bus stops to be allocated at any time. This constraint can be expressed as:
    \begin{align}
        \underline{u}^{(i)}\le \sum_{c_j\in \mathcal{N}}u_{k,c_j}^{(i)} \le \overline{u}^{(i)} \ \forall i \in \mathcal{M}, \ \forall k\in \mathcal{T}\label{eq:max_min_cons}
    \end{align}
    where $\underline{u}^{(i)}$ and $\overline{u}^{(i)}$ are the minimum and maximum number allowable for curb-type $i\in \mathcal{M}$.
    \item Distance between similar allocation types: We would like to spread out similar allocations as much as possible from each other, e.g., we would like to avoid having two bus-stops right next to each other. To achieve this, we can formulate a distance matrix from our spatial data, $A_{\text{m}} \in \mathbb{R}^{N\times N}$, which relates the distance between any two curbs, e.g., the distance between curbs $c_i$ and $c_j$ will be given by $A_{\text{m}}(c_i,c_j)$. If $\mathbf{u_k^{(i)}}$ is the curb allocation vector at time-step $k$ for curb-type $i$, then this constraint can be expressed as:
    \begin{align}
        w_k^{(i)}\le \mathbf{u_k^{(i)}}^\top A_{\text{m}}\mathbf{u_k^{(i)}} \ \forall k\in \mathcal{T}, \forall i \in \mathcal{M}
    \end{align}\label{eq:max_dist_cons}
\end{enumerate}

Based on our objective function and the different constraints that we need to satisfy, the dynamic curb-allocation problem can be formulated as a mixed-integer program (MIP) as follows:
\begin{subequations}\label{eq:MIP_opt}
\begin{align}
    &\max \sum_{k\in \mathcal{T}}\sum_{c_j\in \mathcal{N}}F_k(u_{k,c_j},c_j)+\rho \sum_{k \in \mathcal{T}}\sum_{i \in \mathcal{M}}w_k^{(i)} \label{eq:MIP_a}\\
    &\text{s.t} \  \sum_{c_j \in \mathcal{N}}\sum_{i\in \mathcal{M}}\frac{1}{2}||u_{k+1,c_j}^{(i)}-u_{k,c_j}^{(i)}||_1 \le b, \ \forall k \in \mathcal{T}\backslash T \label{eq:MIP_b}\\
    &\underline{u}^{(i)}\le \sum_{c_j\in \mathcal{N}}u_{k,c_j}^{(i)} \le \overline{u}^{(i)}, \ \forall i \in \mathcal{M}, \ \forall k\in \mathcal{T} \label{eq:MIP_c}\\
    &w_k^{(i)}\le \mathbf{u_k^{(i)}}^\top A_{\text{m}}\mathbf{u_k^{(i)}}, \ \forall k\in \mathcal{T}, \forall i \in \mathcal{M}\label{eq:MIP_d}\\
    &\sum_i u_{k,c_j}^{(i)}=1,\ \forall c_j \in \mathcal{N}, \ \forall k \in \mathcal{T}\label{eq:MIP_e}
\end{align}
\end{subequations}
where the second term in the objective is the regularizer that maximizes the distance between similar allocations, with $\rho$ being the regularizer weight. The constraint in~\eqref{eq:MIP_e} ensures that only one type of curb zoning is allocated to any curb space at each time-step.

Having laid out our MIP, we briefly comment on its relationship to the shared resource multi-agent planning formulation mentioned earlier~\cite{hong2011optimal}. For each curb, i.e., agent, the decision variables are a sequence of integers that represent allocations over time, i.e, a plan. On each individual curb, we have a per-agent constraint of allocating only one type of zoning at a time. Our global objective function decomposes over individual curbs and is the sum of local curb objectives; however, unlike them, we have an additional component (comprising the $w_k^{(i)}$ terms) that acts as a regularizer and depends on inter-agent decisions, i.e., the sum of distances between similar allocations. Finally, our shared inter-agent constraints~\eqref{eq:MIP_b}-\eqref{eq:MIP_d} are more complex than that of the prior work, which imposes a constraint on a sum of linear functions over the agent plans.

\begin{figure}[t]
    \centering
   \subfloat {\includegraphics[width=.48\linewidth]{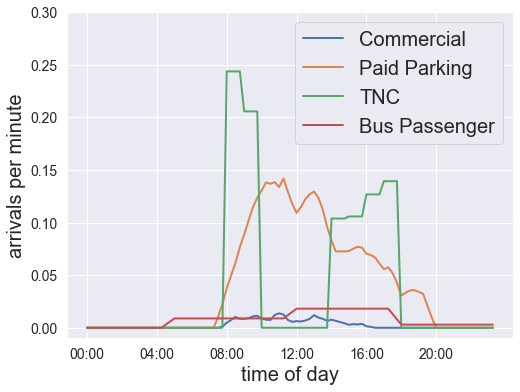}\label{fig:arrival}}
    \hfill
    \subfloat {\includegraphics[width=.48\linewidth]{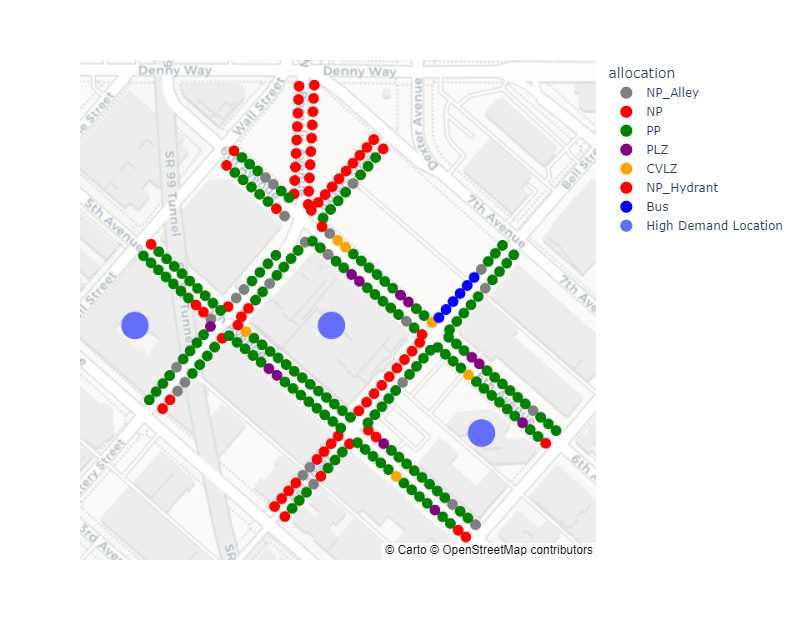}\label{fig:insignia_area}}
    \caption{(a) Average arrival rate of median passengers per vehicle over a typical workday (b) Curb zoning at the time of the data collection study on arrival rate.}
\end{figure}

\subsection{Data and Simulation of MIP}

To showcase the result of the MIP on the curb zoning problem, we consider this problem for a Seattle neighborhood (illustrated in Fig.~\ref{fig:insignia_area}) with 289 curb spaces over a time horizon of 10 hours. We consider three possible curb allocation types to be dynamically rezoned per hour, namely: paid parking (PP), commercial vehicle load zones (CVLZ) and buses/public transportation (Bus). Two primary data points are required to construct an example valuation of individual curb space: arrival rate of vehicles by type (an implicit measure of demand), and the parking space's distance from their intended destination (where closer is more desirable \cite{de2018parking}). For each curb space in the area of interest there are many ways to \emph{value} the real estate: for this example we compute a net revenue to the municipality, both internalized (e.g., price to park, value of time spent walking to destination), and externalized (e.g., cost of carbon, cost of congestion). All data, relevant sources, and code are available at \url{https://github.com/cpatdowling/dynamiczoning}.

\begin{figure}[t]
    \centering
   \subfloat {\includegraphics[width=.48\linewidth]{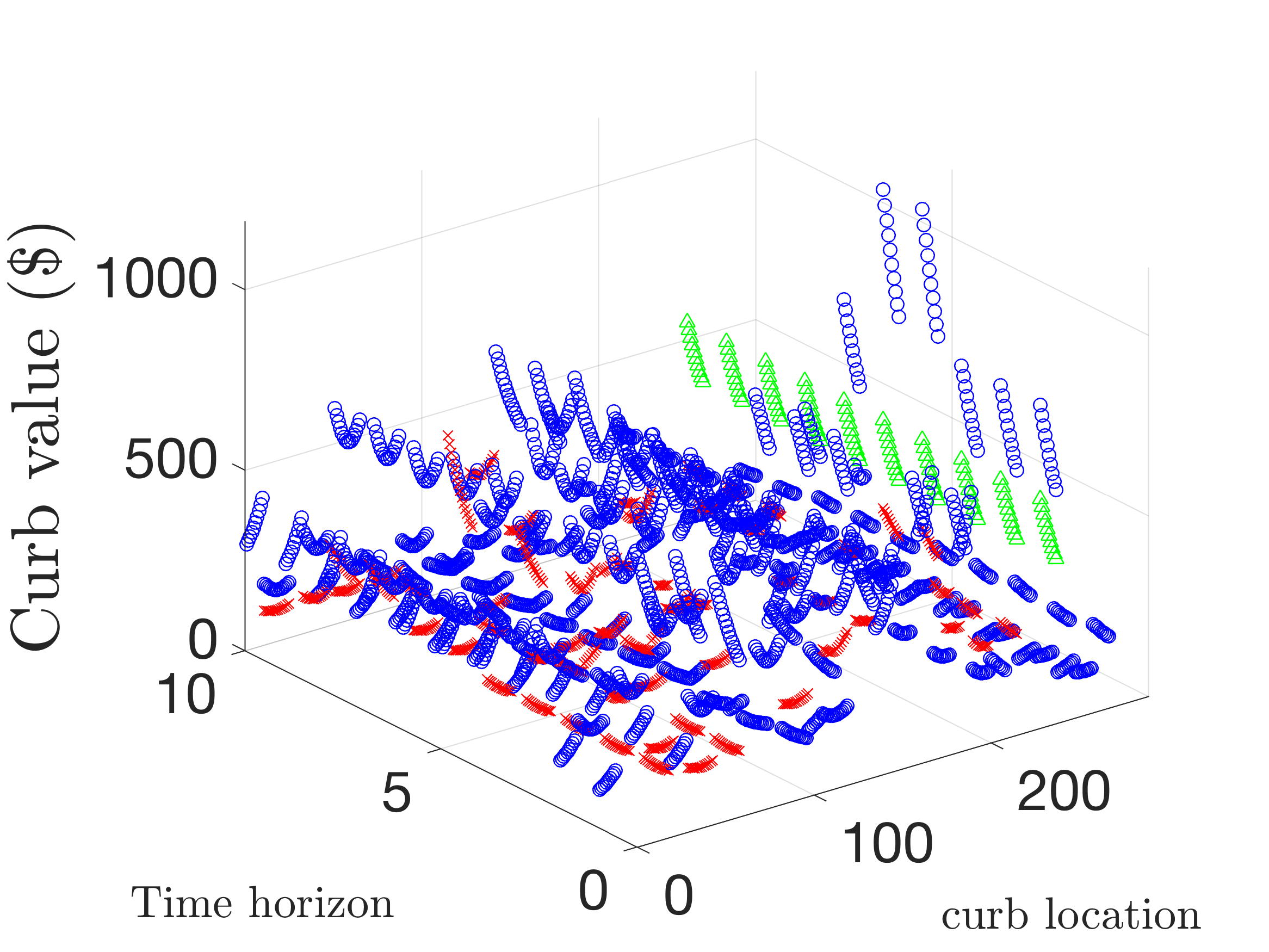}\label{fig:curb_value}}
    \hfill
    \subfloat {\includegraphics[width=.48\linewidth]{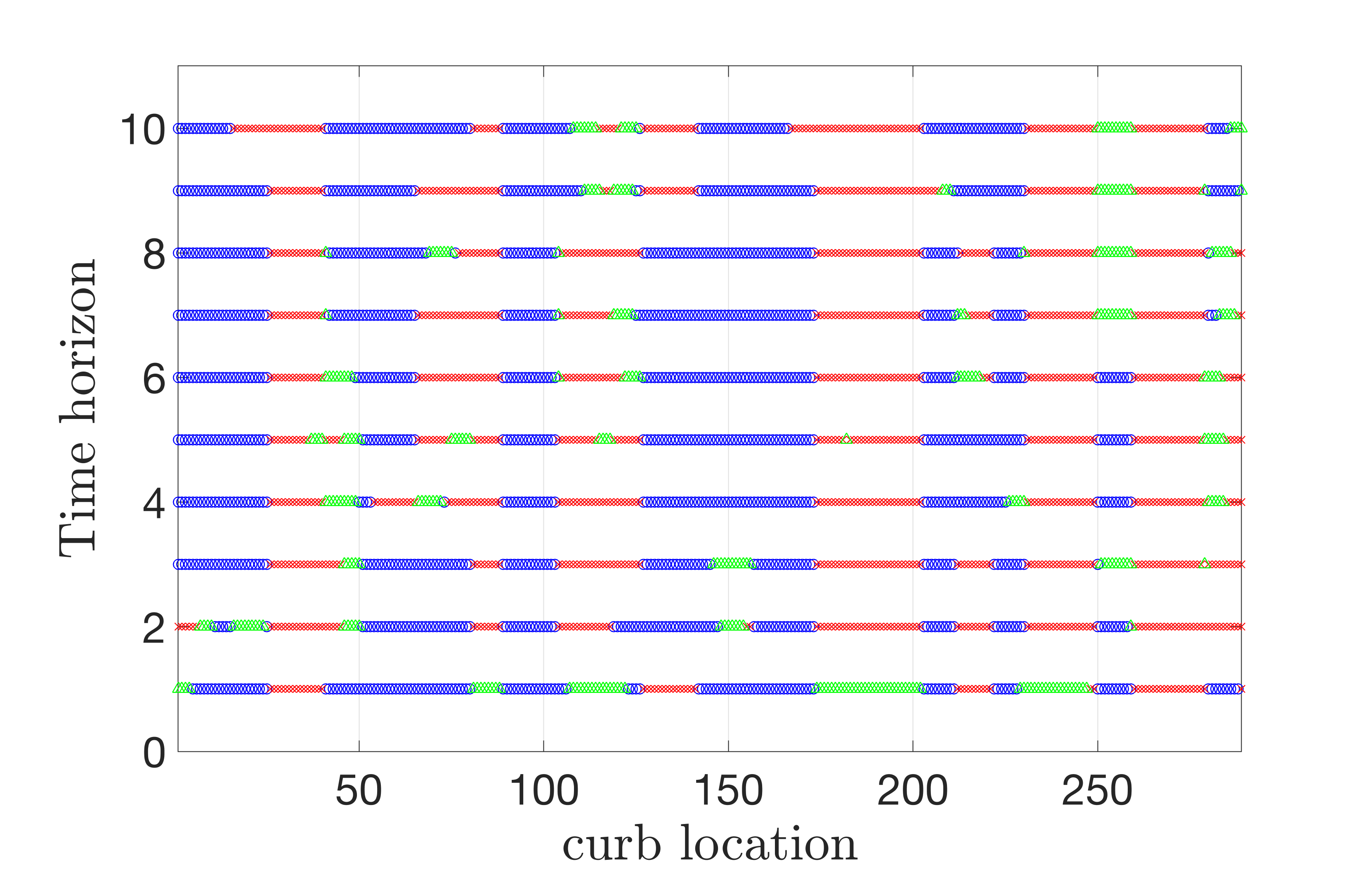}\label{fig:MIP_allocation}}
    \caption{(a) Value of curb spaces over time and curb locations for three different types of allocations: blue-paid parking, red-commercial parking and green-buses/public transportation (b) Allocation of curbs at different locations and at different times for different zoning purposes based on solution of the MIP.}
\end{figure}

The key data source in approximating time-dependent demand for space in our example is the arrival rate of passengers (in the case of bus transit, paid parking, and passenger pickup/drop-off) and vehicles (in the case of commercial vehicles). To measure these arrival rates, we combine temporally coinciding paid parking transactions, \cite{SDOT2019, dowling2019modeling}, King County Metro study data, and manual observation studies \cite{giron2019commercial}, combined in Fig.~\ref{fig:arrival}. Then for each individual parking space location for each zone type, we compute a valuation based on historically expected demand and distance to hypothetical points of interest \cite{nie2021public}. The value for each allocation at each of the curb spaces over the time horizon is shown in Fig.~\ref{fig:curb_value}. Based on this curb value, the MIP optimization problem~\eqref{eq:MIP_opt} is solved to achieve the optimal curb zoning. The results of this optimization solution are depicted in Fig.~\ref{fig:MIP_allocation} that shows the curb zoning over the curb spaces and over the time horizon and satisfying the given network and operational constraints.

To be clear, while we have attempted to construct an example objective function with real data, we would strongly encourage practitioners to consider more principled research focused on measuring the value of various zoning types \cite{van2011real, shoup2021pricing}. This work is focused on applying such a measurement to optimally allocate or rezone curb space over time, and showing that this is a dynamic programming problem.


\subsection{Drawbacks of MIP}
The above MIP can be solved with many commercially available MIP solvers. However, due to the large scale nature of the problem, it becomes difficult to scale up when dealing with large number of curb spaces in a city as the solution time for integer programs grows exponentially~\cite{lodi2010mixed}. In order to deal with this scalability issue, we instead rely on utilizing two different techniques to solve MIPs. In the first approach, we utilize certain approximate dynamic programming (ADP) solution methods. The drawback of ADP is that it provides no guarantees on global optimality. In the second approach, we utilize a Dantzig-Wolfe decomposition method to solve this MIP. In the next sections, we will compare different ADP methods to Dantzig-Wolfe on their computation time and the optimality gap with the MIP solution.

\section{Approximate dynamic programming}\label{sec:ADP}


The mixed-integer programming based approach is not scalable as the problem size grows exponentially with time horizon. To address this we utilize dynamic programming techniques to solve the curb allocation problem. Specifically we aim to utilize recent advances in approximate dynamic programming (ADP) to achieve scalable yet sufficiently accurate solutions to the curb allocation problem.
Dynamic programming breaks a multi-period planning problem into simpler sub-problems at different points in time. This optimization problem is then solved in a recursive manner by using the Bellman equation.
Let the state of the system at time-step $k$ be $x_k$, with $x_0$ being the initial state. Then, an infinite-horizon decision problem takes the following form:
\begin{align}
    V_0(x_0)=\max_{u_k}\sum_{k=0}^{\infty}F_k(u_k)  : \ x_{k+1}=G(x_k,u_k)\label{eq:DP_1}
\end{align}
where $V_0(x_0)$ is the optimal value of the objective function for initial state $x_0$, $u_k$ is the decision variable at time-step $k$, $F_k(u_k)$ is the objective function at time-step $k$ and  $G(x_k,u_k)$ represents the system state-change dynamics. One can write the above problem in a recursive manner as:
\begin{align}
    V_0(x_0)=\max_{u_0}\{F_0(u_0)+V_1(x_1)\}  : \ x_1=G(x_0,u_0) \label{eq:DP_2}
\end{align}
which can be written in the generalized form as:
\begin{equation}
\begin{aligned}
    &V_k(x_k)=\max_{u_k}\{F_k(u_k)+V_{k+1}(x_{k+1})\}\\
    &x_{k+1}=G(x_k,u_k) \label{eq:DP_3}
\end{aligned}
\end{equation}
Since in this case, the objective function satisfies a separable structure~\cite{jones2021generalization} and $u_k$ does not depend upon $u_{k+1}$, the Bellman equation can be written without the state variables as:

\begin{align}
    V_0=\sum_{k=0}^{\infty}\max_{u_k}F_k(u_k) \label{eq:DP_sep}
\end{align}
This special separable structure of the objective function allows us to compute the optimal value function at each time-step independently. However, we still need to consider the inter-temporal coupling constraints in~\eqref{eq:MIP_opt}. We use these simplifications to develop the ADP formulation shown next.

\begin{algorithm}[t]
\caption{\label{alg1} Adaptive dynamic programming (ADP4) based curb allocation algorithm}
\SetAlgoLined
\KwResult{optimal allocation $u_{k,c_j}^{(i)}$}
\textbf{Input:} $F_k, b, A_{\text{m}},\underline{u}^{(i)}, \overline{u}^{(i)}$\\
 Initialize random $u_{k,c_j}^{(i)}$\\
 \For{$p=1:N_{OL}$}{
 \For{$k=1:T$}{
 \For{$q=1:N_{IL}$}{
Choose 10 random curb locations\\
Modify $u_{k,c_j}$ for these locations based on maximum $F_k$\\
Calculate objective value for this allocation: $V_k(q)=F_k(q)+\rho w_{k}(q)$\\
    \If{(2) or (3) not satisfied}{
    $V_k(q)=-\infty$
}
}
 $u_{k,c_j}=\argmax (V_k(q)$) \\
}
Save optimal allocation for comparison in next iteration\\
}
\end{algorithm}

In the ADP formulation, we utilize Monte Carlo methods to sample the feasible space and iteratively converge towards an optimal solution~\cite{powell2007approximate}. Because of the nature of curb demand distribution, however, the optimal curb allocation solution tends to consist of blocks of similar curb allocations together in groups. We can also ascertain this from the clusters in the solution to the MIP illustrated in Fig.~\ref{fig:MIP_allocation}. 

\begin{figure}[t]
    \centering
    \includegraphics[width=\linewidth]{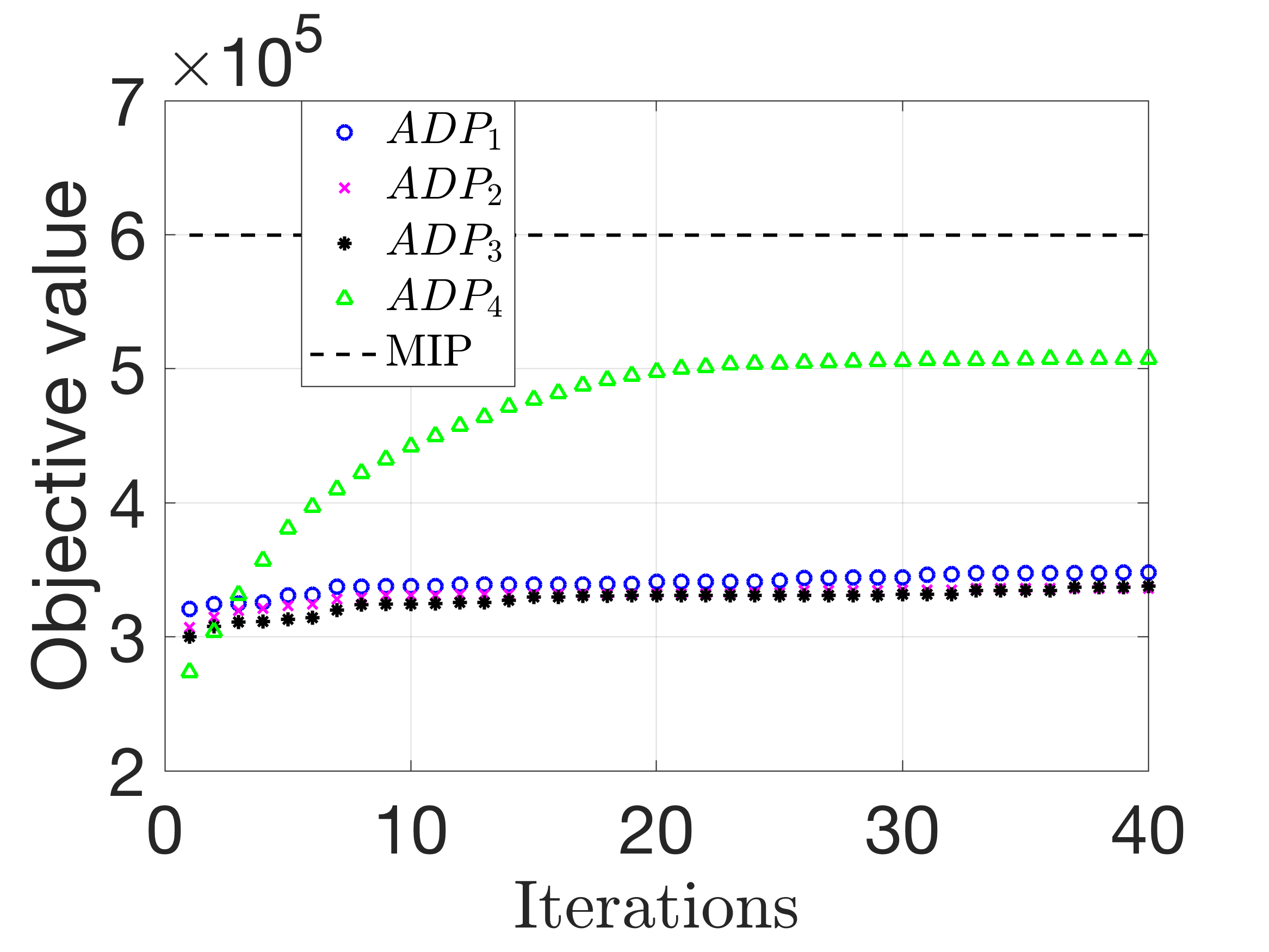}
    \caption{Comparison of different ADP techniques with the MIP solution.\\
     $ADP_1$: Randomly generate the total number of curbs allocated to paid parking, commercial vehicles and buses and update that many locations, calculate the value function, save the maximum value function allocation and repeat.\\
    $ADP_2$: Randomly choose certain curbs and randomly update their allocations (between paid parking, commercial vehicles and buses), calculate value function, save the maximum value function allocation and repeat.\\
    $ADP_3$: Same as $ADP_2$ but every often choose a completely random allocation, begin again and repeat $ADP_2$.\\
    $ADP_4$: Randomly choose certain number of curbs and update their allocation according to local maximum curb value function, then calculate total value function and save the maximum. At batched intervals, choose a completely random allocation and repeat again.}
    \label{fig:ADP_comp}
\end{figure}

Purely random Monte Carlo methods are unable to sample the feasible space in such a clustered manner. Thus, we formulate a hybrid Monte Carlo algorithm which considers both the global and the local nature of the cost function when allocating curb zones. The algorithm shown in Algorithm~\ref{alg1} consists of two loops. In the outer loop we randomly sample the feasible space and converge towards an optimal solution. In the inner loop we consider the local cost of each curb space and allocate them according to the localized costs. This mechanism can be thought of as a combination of global exploration and local exploitation, where we locally allocate curb spaces based on their cost functions and then explore random allocations to compare with. Allocations that maximize the objective along with satisfying all the network constraints are stored as optimal. Then the process repeats and new solutions are compared with the stored optimal value and the stored optimal is updated if a better solution is obtained. The process is repeated until we converge to a solution. Algorithm~\ref{alg1} below shows the architecture of the proposed algorithm. Then Fig.~\ref{fig:ADP_comp} compares the results of our algorithm (ADP4) with those of simple Monte Carlo based ADPs (ADP1, ADP2, ADP3) and the globally optimal solution to the full MIP.

\section{Dantzig-Wolfe decomposition}\label{sec:DW}

A fundamental approach for solving linear programs (LPs), Dantzig-Wolfe decomposes the overall problem into a \emph{master} problem that imposes inter-agent shared resource constraints, along with one \emph{subproblem} for each agent~\cite{bertsimas1997introduction}.
The master problem has far fewer constraints than the original, as it ignores per-agent constraints. But its domain is the set of all feasible multi-agent plans and is thus exponentially large in the number of agents. Therefore, Dantzig-Wolfe starts off the master problem with a restricted solution set that it expands iteratively with \emph{delayed column generation} by solving subproblems, terminating only when a subproblem can no longer contribute a better feasible solution to the master problem.

To apply Dantzig-Wolfe decomposition to our specific curb-allocation setting, we exploit its special structure in a manner similar to the prior work on multi-agent planning~\cite{hong2011optimal}. We relax the integer variables and solve an LP as the master problem and obtain the solution which then solves the curb-allocation problem for each curb space in parallel, considering only the local curb constraints along with the broadcasted duals. At each of the sub-problems, we calculate the reduced cost, and if it is negative, we add that allocation to the list of possible feasible solutions for when we next re-solve the master problem. This process is repeated in an iterative fashion until the solution converges.

The relaxed LP is shown below:

\begin{subequations}\label{eq:DW_LP}
\begin{align}
    &\max \sum_{k\in \mathcal{T}}\sum_{c_j\in \mathcal{N}}F_k(u^{(m)}_{k,c_j},c_j) \label{eq:DW_LP_obj}\\
    \begin{split}
     &\text{s.t} \   \sum_{c_j \in \mathcal{N}}\sum_{i\in \mathcal{M}}\frac{1}{2}||u_{k+1,c_j}^{(i)}-u_{k,c_j}^{(i)}||_1 \le b,\forall k \in \mathcal{T}\backslash T : \nu_k\label{eq:diff_constraint}
     \end{split}\\
     \begin{split}
    &\underline{u}^{(i)}\le \sum_{c_j\in \mathcal{N}}u_{k,c_j}^{(m)(i)} \le \overline{u}^{(i)},\forall i \in \mathcal{M}, k\in \mathcal{T} : [\underline{\lambda}_k^{(i)},\overline{\lambda}_k^{(i)}] \label{eq:lagrange_bounds}
    \end{split}\\
    &u^{(m)(i)}_{k,c_j}=\sum_{m \in \mathcal{Q}}z_mU^{(m)}\label{eq:Um_constraint}\\
    &\sum_{m\in \mathcal{Q}}z_m=1 \ : \ \mu \label{eq:zsum_constraint}\\
    &z_m\in \mathcal{Z}, \ \forall m \in \mathcal{Q} \label{eq:MIP_constraint_LP}
\end{align}
\end{subequations}\label{eq:LP_opt}

where $\mathcal{Q}$ is the set of all feasible solutions $(m)$ to the original MIP, $U^{m}$ is the curb allocation for  the feasible solutions and $[\underline{\lambda}_k^{(i)},\overline{\lambda}_k^{(i)}]$ are the lagrange duals associated with the constraints in~\eqref{eq:lagrange_bounds}, $\mu$ is the lagrange dual associated with the constraint in~\eqref{eq:zsum_constraint} and $\nu_k$ is the dual associated with the constraint in~\eqref{eq:diff_constraint}. \textcolor{black}{Here, \eqref{eq:Um_constraint} chooses the binary weights associated with the curb allocation of feasible solutions $U^{m}$ and \eqref{eq:zsum_constraint} ensures that only one set of feasible solutions is chosen.} The above problem can be turned into an LP by relaxing the integer constraint in~\eqref{eq:MIP_constraint_LP} to the following linear constraint:

\begin{align}
    z_m \in [0,1], \ \forall m\in \mathcal{Q}\label{eq:relax_cons}
\end{align}

Utilizing the Lagrange duals obtained from solving the above LP, these values are broadcasted to solve a MIP locally and in parallel for each curb location $c_j$ as follows:

\begin{subequations}
\begin{align}
\begin{split}
    &\max \sum_{k\in \mathcal{T}}F_k(u_{k,c_j},c_j)-\sum_{i \in \mathcal{M}}\sum_{k \in \mathcal{T}}\overline{\lambda}_k^{(i)}u^{(i)}_{k,c_j}-\\
   & \sum_{i \in \mathcal{M}}\sum_{k \in \mathcal{T}}\underline{\lambda}_k^{(i)}u^{(i)}_{k,c_j}-\sum_{k\in \mathcal{T}\backslash T}\sum_{i\in \mathcal{M}}\nu_ku^{(i)}_{k+1,c_j}
   \end{split}\\
    &\text{s.t} \ \sum_{i \in \mathcal{M}}u_{k,c_j}^{(i)}=1, \ \forall k \in \mathcal{T}\\
   & u_{k,c_j}^{(i)} \in \mathcal{Z}
\end{align}
\end{subequations}\label{eq:MIP_local}

The optimal solution obtained from solving the above MIP is used to calculate the reduced cost shown below:

\begin{align}
\begin{split}
    &\sum_{k\in \mathcal{T}}F_k(u^*_{k,c_j},c_j)-\sum_{i \in \mathcal{M}}\sum_{k \in \mathcal{T}}\overline{\lambda}_k^{(i)}u^{*(i)}_{k,c_j}-\\
    &\sum_{i \in \mathcal{M}}\sum_{k \in \mathcal{T}}\underline{\lambda}_k^{(i)}u^{*(i)}_{k,c_j}-\sum_{k\in \mathcal{T}\backslash T}\sum_{i\in \mathcal{M}}\nu_ku^{*(i)}_{k+1,c_j}-\mu
    \end{split}
\end{align}

If the reduced cost is less than zero, then the obtained optimal solution is added to the list of feasible solutions $U^m$ for the LP and the LP is resolved with the addition of the new feasible solution. This process is repeated until convergence.

\begin{figure}
    \centering
    \includegraphics[width=.5\textwidth]{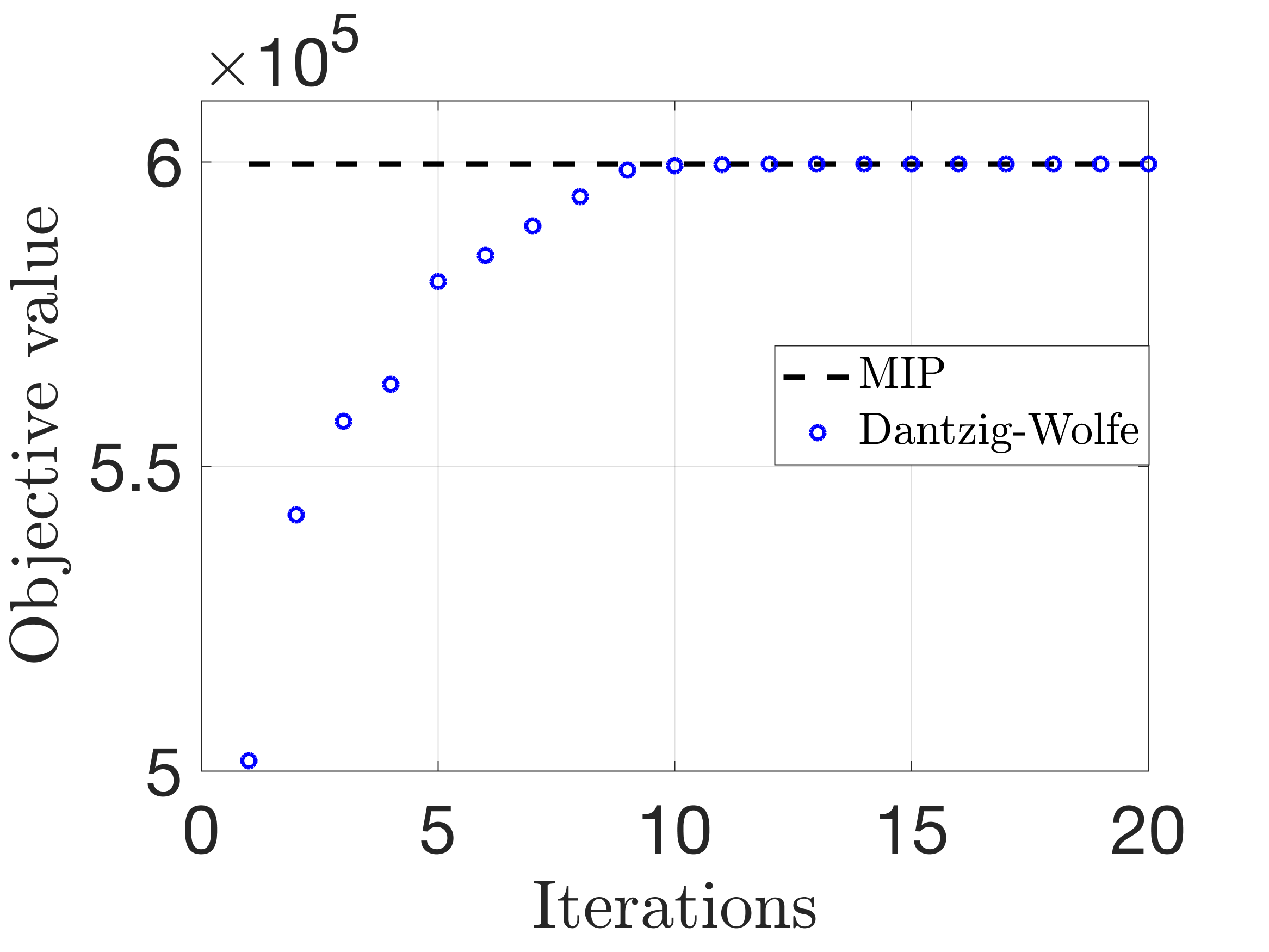}
    \caption{Comparison of objective value obtained using Dantzig-Wolfe decomposition with the MIP solution over iterations.}
    \label{fig:DW_comp_Seattle}
\end{figure}

Fig.~\ref{fig:DW_comp_Seattle} compares the simulation results using the Dantzig-Wolfe approach to the MIP solution, for the Seattle insignia neighborhood with a total of 289 curb spaces over 10 time-steps. From these results, we are able to demonstrate the efficacy of Dantzig-Wolfe in converging towards the optimal solution obtained from the MIP. Of course, we are primarily motivated to use Dantzig-Wolfe decomposition due to its scalability, as it enables us to solve an LP combined with smaller MIPs in a parallel manner. This circumvents the scalability issue encountered in solving large-scale MIPs.

\begin{figure}
    \centering
    \includegraphics[width=.5\textwidth]{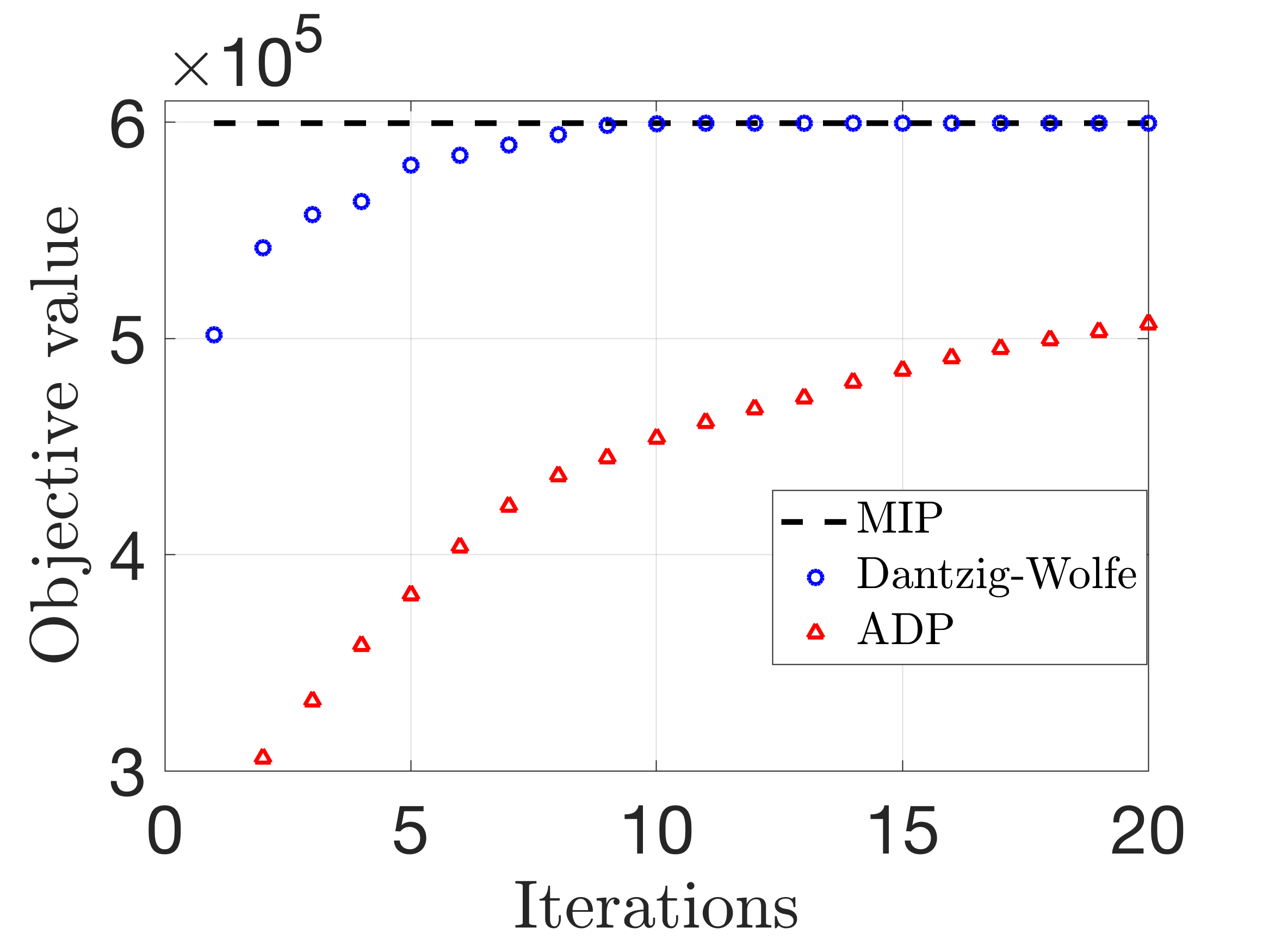}
    \caption{Comparison of objective value obtained using Dantzig-Wolfe decomposition with the ADP4 and MIP solution over iterations.}
    \label{fig:DW_ADP_comp_Seattle}
\end{figure}

\begin{figure}
    \centering
    \includegraphics[width=.5\textwidth]{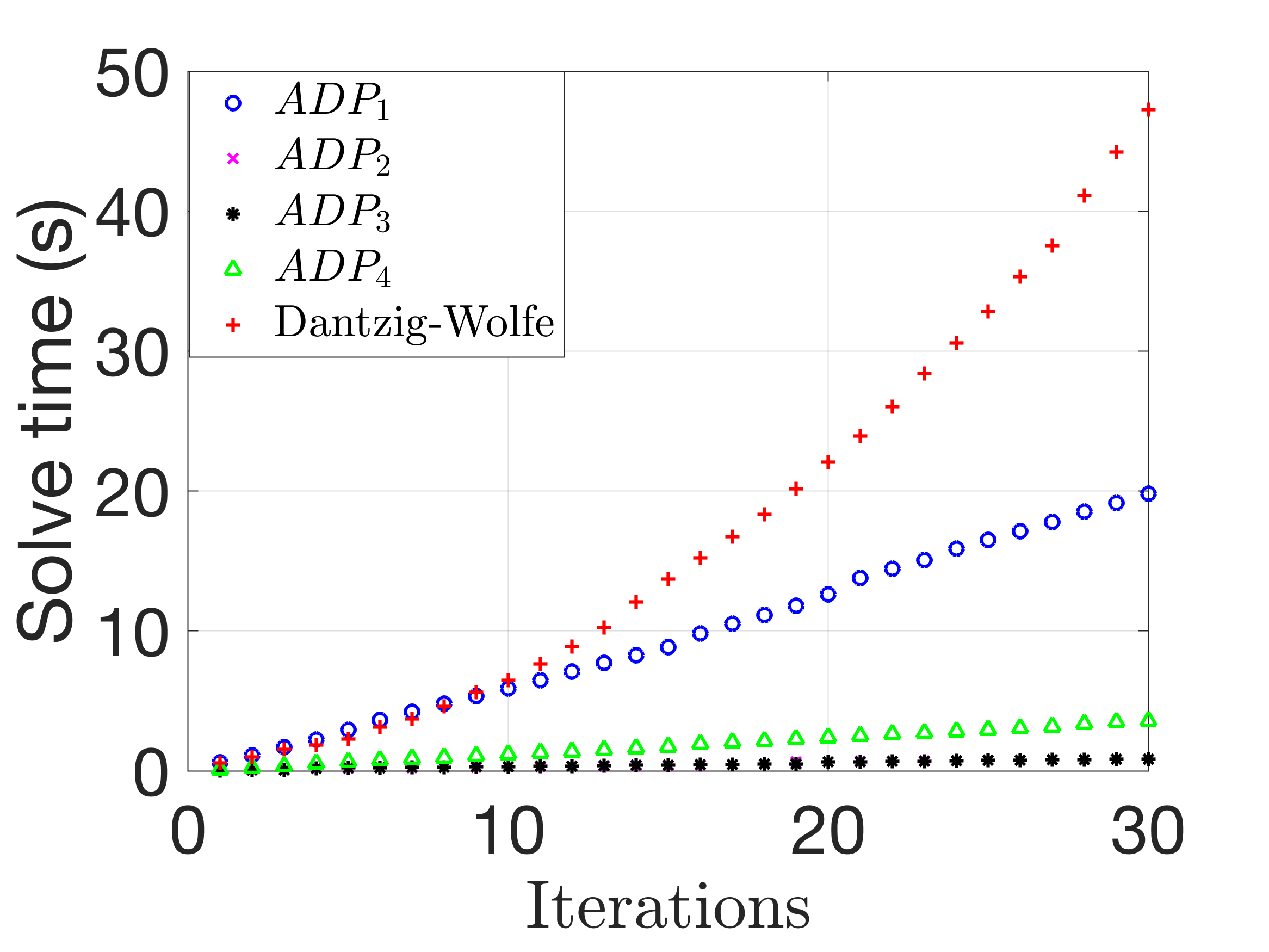}
    \caption{Comparison of cumulative solve time over iterations between Dantzig-Wolfe decomposition method with the ADP methods.}
    \label{fig:DW_ADP_comp_Seattle_time}
\end{figure}


\section{Discussion}\label{sec:Discussion}

A critical assumption of this work is the existence of a curb value function that a municipality can reasonably measure or approximate, to then optimize over. In this work, we have maximized an objective that encodes the net revenue of the curb to all users relative to the distance from hypothetical points of interest across three types: passenger parking, commercial loading, and transit riders (with constraints on the minimum and maximum number of each zone type, maximum number of rezonings over the time horizon, and a penalty on solution that deviate from a desired distance between bus stops). We see two pathways forward. First, a wide body of econometric work seeks to measure these valuations directly, particularly with respect to the externalized costs of a given zoning type~\cite{van2011real, de2018parking, van2018marginal, shoup2021pricing}. Second, one could approximate an implicit valuation by measuring the changing response in demand for a space when a space is rezoned. Such an approach would build upon inverse reinforcement learning \cite{ng2000algorithms} and could incorporate domain and institutional knowledge through Bayesian priors~\cite{ramachandran2007bayesian}.

The latter introduces an interesting problem that may be considered by more dynamic management of future curb zoning: as curb space is rezoned, demand at that space and for spaces in neighboring areas may change as a result. Resolving this form of control-dependent demand might be addressed through performative prediction methods in the context of our proposed dynamic programming formulation \cite{narang2022multiplayer}.

\section{Conclusions and Future work}\label{sec:Conclusions}

In this paper we have shown that, given some principled measure of curb valuation, dynamically rezoning curb space is a dynamic programming problem. We also show that we can easily introduce policy-driven constraints, such as upper and lower bounds on the number of a given zoning type, or a desired distance between any two spaces zoned with the same type (e.g., bus stops). For even modestly large neighborhoods, however, solving this dynamic programming problem becomes computationally intractable \emph{particularly} if these policy-driven constraints are coupled in time, such as limiting the total number of rezonings over a finite time period. We show that this can be addressed by sampling the solution space through the use of approximate dynamic programming, or by decomposing the problem via Dantzig-Wolfe decomposition.

\textcolor{black}{Future work could model other types of policy-related time constraints in the curbside zoning problem. It could also demonstrate our approach on case studies with more zone types and longer planning
horizons. Finally, another important direction for future work will be to conduct a sensitivity analysis of optimal curb-zoning
with respect to various system parameters, in order to yields\ helpful insight on real-world implementation.}

\section*{Acknowledgements}
Pacific Northwest National Laboratory is operated by Battelle Memorial Institute for the U.S. Department of Energy under Contract No. DE-AC05-76RL01830. This work was supported by the U.S. Department of Energy Vehicle Technologies Office.






\bibliographystyle{IEEEtran}
\bibliography{sample}

\begin{thebibliography}{10}
\providecommand{\url}[1]{#1}
\csname url@samestyle\endcsname
\providecommand{\newblock}{\relax}
\providecommand{\bibinfo}[2]{#2}
\providecommand{\BIBentrySTDinterwordspacing}{\spaceskip=0pt\relax}
\providecommand{\BIBentryALTinterwordstretchfactor}{4}
\providecommand{\BIBentryALTinterwordspacing}{\spaceskip=\fontdimen2\font plus
\BIBentryALTinterwordstretchfactor\fontdimen3\font minus
  \fontdimen4\font\relax}
\providecommand{\BIBforeignlanguage}[2]{{%
\expandafter\ifx\csname l@#1\endcsname\relax
\typeout{** WARNING: IEEEtran.bst: No hyphenation pattern has been}%
\typeout{** loaded for the language `#1'. Using the pattern for}%
\typeout{** the default language instead.}%
\else
\language=\csname l@#1\endcsname
\fi
#2}}
\providecommand{\BIBdecl}{\relax}
\BIBdecl

\bibitem{abel2021curbside}
S.~Abel, M.~Ballard, S.~Davis, M.~Mitman, K.~Stangl, D.~Wasserman
  \emph{et~al.}, ``Curbside inventory report,'' United States. Federal Highway
  Administration, Tech. Rep., 2021.

\bibitem{blakeley2013time}
M.~Blakeley and N.~Gray, ``Time-lapse cameras measure street parking demand,''
  \emph{Institute of Transportation Engineers. ITE Journal}, vol.~83, no.~9,
  p.~36, 2013.

\bibitem{ranjbari2021testing}
A.~Ranjbari, J.~Luis Machado-Le{\'o}n, G.~Dalla~Chiara, D.~MacKenzie, and
  A.~Goodchild, ``Testing curbside management strategies to mitigate the
  impacts of ridesourcing services on traffic,'' \emph{Transportation Research
  Record}, vol. 2675, no.~2, pp. 219--232, 2021.

\bibitem{saha2019project}
M.~Saha, M.~Saugstad, H.~T. Maddali, A.~Zeng, R.~Holland, S.~Bower, A.~Dash,
  S.~Chen, A.~Li, K.~Hara \emph{et~al.}, ``Project sidewalk: A web-based
  crowdsourcing tool for collecting sidewalk accessibility data at scale,'' in
  \emph{Proceedings of the 2019 CHI Conference on Human Factors in Computing
  Systems}, 2019, pp. 1--14.

\bibitem{fiez2018data}
T.~Fiez, L.~J. Ratliff, C.~Dowling, and B.~Zhang, ``Data driven spatio-temporal
  modeling of parking demand,'' in \emph{2018 Annual American Control
  Conference (ACC)}.\hskip 1em plus 0.5em minus 0.4em\relax IEEE, 2018, pp.
  2757--2762.

\bibitem{roe2017curb}
M.~Roe and C.~Toocheck, ``Curb appeal: Curbside management strategies for
  improving transit reliability,'' \emph{National Association of City
  Transportation Officials}, 2017.

\bibitem{harris2017mesoscopic}
T.~M. Harris, M.~Nourinejad, and M.~J. Roorda, ``A mesoscopic simulation model
  for airport curbside management,'' \emph{Journal of Advanced Transportation},
  vol. 2017, 2017.

\bibitem{giron2019commercial}
G.~d.~C. Gir{\'o}n-Valderrama, J.~L. Machado-Leon, and A.~Goodchild,
  ``Commercial vehicle parking in downtown seattle: insights on the battle for
  the curb,'' \emph{Transportation Research Record}, vol. 2673, no.~10, pp.
  770--780, 2019.

\bibitem{yu2021management}
M.~Yu and A.~Bayram, ``Management of the curb space allocation in urban
  transportation system,'' \emph{International Transactions in Operational
  Research}, vol.~28, no.~5, pp. 2414--2439, 2021.

\bibitem{schwaeppe2019equal}
H.~Schwaeppe, M.~Nobis, and C.~M{\"u}ller, ``On the equal substitution of milp
  unit commitment subproblems with dynamic programming,'' in \emph{2019 16th
  International Conference on the European Energy Market (EEM)}.\hskip 1em plus
  0.5em minus 0.4em\relax IEEE, 2019, pp. 1--5.

\bibitem{lewis2009reinforcement}
F.~L. Lewis and D.~Vrabie, ``Reinforcement learning and adaptive dynamic
  programming for feedback control,'' \emph{IEEE circuits and systems
  magazine}, vol.~9, no.~3, pp. 32--50, 2009.

\bibitem{powell2007approximate}
W.~B. Powell, \emph{Approximate Dynamic Programming: Solving the curses of
  dimensionality}.\hskip 1em plus 0.5em minus 0.4em\relax John Wiley \& Sons,
  2007, vol. 703.

\bibitem{barnhart1998branch}
C.~Barnhart, E.~L. Johnson, G.~L. Nemhauser, M.~W. Savelsbergh, and P.~H.
  Vance, ``Branch-and-price: Column generation for solving huge integer
  programs,'' \emph{Operations research}, vol.~46, no.~3, pp. 316--329, 1998.

\bibitem{galati2010decomposition}
M.~Galati, \emph{Decomposition methods for integer linear programming}.\hskip
  1em plus 0.5em minus 0.4em\relax Lehigh University, 2010.

\bibitem{andrianesis2021computation}
P.~Andrianesis, D.~J. Bertsimas, M.~Caramanis, and W.~Hogan, ``Computation of
  convex hull prices in electricity markets with non-convexities using
  dantzig-wolfe decomposition,'' \emph{IEEE Transactions on Power Systems},
  2021.

\bibitem{guide2019introducing}
P.~GUIDE, ``Introducing ite’s new curbside management practitioners guide,''
  \emph{Institute of Transportation Engineers}, 2019.

\bibitem{ostermeijer2021citywide}
F.~Ostermeijer, H.~Koster, L.~Lunes, and J.~N. van Ommeren, ``Citywide parking
  policy and traffic: Evidence from amsterdam,'' \emph{Journal of Urban
  Economics}, 2021.

\bibitem{inci2015review}
E.~Inci, ``A review of the economics of parking,'' \emph{Economics of
  Transportation}, vol.~4, no. 1-2, pp. 50--63, 2015.

\bibitem{kong2018iot}
X.~T. Kong, S.~X. Xu, M.~Cheng, and G.~Q. Huang, ``Iot-enabled parking space
  sharing and allocation mechanisms,'' \emph{IEEE Transactions on Automation
  Science and Engineering}, vol.~15, no.~4, pp. 1654--1664, 2018.

\bibitem{nakazato2019parking}
T.~Nakazato and T.~Namerikawa, ``Parking lot allocation based on matching
  theory using prediction-based optimal vehicle routing,'' in \emph{2019 19th
  International Conference on Control, Automation and Systems (ICCAS)}.\hskip
  1em plus 0.5em minus 0.4em\relax IEEE, 2019, pp. 1004--1009.

\bibitem{torreno2017cooperative}
A.~Torre{\~n}o, E.~Onaindia, A.~Komenda, and M.~{\v{S}}tolba, ``Cooperative
  multi-agent planning: A survey,'' \emph{ACM Computing Surveys (CSUR)},
  vol.~50, no.~6, pp. 1--32, 2017.

\bibitem{hong2011optimal}
S.~A. Hong and G.~Gordon, ``Optimal distributed market-based planning for
  multi-agent systems with shared resources,'' in \emph{Proceedings of the
  Fourteenth International Conference on Artificial Intelligence and
  Statistics}.\hskip 1em plus 0.5em minus 0.4em\relax JMLR Workshop and
  Conference Proceedings, 2011, pp. 351--360.

\bibitem{de2018parking}
D.~de~Vos and J.~Van~Ommeren, ``Parking occupancy and external walking costs in
  residential parking areas,'' \emph{Journal of Transport Economics and Policy
  (JTEP)}, vol.~52, no.~3, pp. 221--238, 2018.

\bibitem{SDOT2019}
SDOT, ``Paid occupancy [dataset],''
  \url{https://data.seattle.gov/Transportation/2020-Paid-Parking-Occupancy-Year-to-date-/wtpb-jp8d},
  2019.

\bibitem{dowling2019modeling}
C.~P. Dowling, L.~J. Ratliff, and B.~Zhang, ``Modeling curbside parking as a
  network of finite capacity queues,'' \emph{IEEE Transactions on Intelligent
  Transportation Systems}, vol.~21, no.~3, pp. 1011--1022, 2019.

\bibitem{nie2021public}
Y.~Nie, W.~Yang, Z.~Chen, N.~Lu, L.~Huang, and H.~Huang, ``Public curb parking
  demand estimation with poi distribution,'' \emph{IEEE Transactions on
  Intelligent Transportation Systems}, 2021.

\bibitem{van2011real}
J.~Van~Ommeren, D.~Wentink, and J.~Dekkers, ``The real price of parking
  policy,'' \emph{Journal of Urban Economics}, vol.~70, no.~1, pp. 25--31,
  2011.

\bibitem{shoup2021pricing}
D.~Shoup, ``Pricing curb parking,'' \emph{Transportation Research Part A:
  Policy and Practice}, vol. 154, pp. 399--412, 2021.

\bibitem{lodi2010mixed}
A.~Lodi, ``Mixed integer programming computation,'' in \emph{50 years of
  integer programming 1958-2008}.\hskip 1em plus 0.5em minus 0.4em\relax
  Springer, 2010, pp. 619--645.

\bibitem{jones2021generalization}
M.~Jones and M.~M. Peet, ``A generalization of bellman’s equation with
  application to path planning, obstacle avoidance and invariant set
  estimation,'' \emph{Automatica}, vol. 127, p. 109510, 2021.

\bibitem{bertsimas1997introduction}
D.~Bertsimas and J.~N. Tsitsiklis, \emph{Introduction to linear
  optimization}.\hskip 1em plus 0.5em minus 0.4em\relax Athena Scientific
  Belmont, MA, 1997, vol.~6.

\bibitem{van2018marginal}
J.~van Ommeren and M.~McIvor, ``The marginal external cost of street parking,
  optimal pricing and supply: evidence from melbourne,'' in \emph{ITEA Annual
  Conference, held in Hong Kong}, 2018, p.~28.

\bibitem{ng2000algorithms}
A.~Y. Ng, S.~J. Russell \emph{et~al.}, ``Algorithms for inverse reinforcement
  learning.'' in \emph{Icml}, vol.~1, 2000, p.~2.

\bibitem{ramachandran2007bayesian}
D.~Ramachandran and E.~Amir, ``Bayesian inverse reinforcement learning.'' in
  \emph{IJCAI}, vol.~7, 2007, pp. 2586--2591.

\bibitem{narang2022multiplayer}
A.~Narang, E.~Faulkner, D.~Drusvyatskiy, M.~Fazel, and L.~J. Ratliff,
  ``Multiplayer performative prediction: Learning in decision-dependent
  games,'' \emph{arXiv preprint arXiv:2201.03398}, 2022.

\end{thebibliography}
\end{document}